\documentclass[a4paper,12pt]{article}

\usepackage{amsthm}
\usepackage{amsmath,amssymb,latexsym,amsfonts,mathrsfs}
\usepackage[dvipdfmx]{graphicx}

\usepackage{color}

\usepackage{multicol}

\usepackage{fancyhdr}
\usepackage[top=30truemm, bottom=30truemm, left=25truemm, right=25truemm]{geometry}

\newcommand{\nn}{\nonumber}

\newcommand{\SCR}[1]{{\mathscr #1}}

\newcommand{\CAL}[1]{{\cal #1}
}

\newcommand{\J}[1]{\left\langle #1 \right\rangle}
\newcommand{\D}[1]{{\mathscr D}( #1 )}

\theoremstyle{definition}
\newtheorem{Thm}{{\bf Theorem}}[section]

\newtheorem{Lem}[Thm]{{\bf Lemma}}
\newtheorem{Prop}[Thm]{{\bf Proposition}}

\newtheorem{Rem}[Thm]{{\bf Remark}}

\newcounter{Exami}

\newcommand{\Proof}[2][Proof]{
\begin{proof}[{\bf #1}]
#2
\end{proof}
}

\begin{document}

\begin{flushleft}
{\bf \Large $L^2$-stableness for solution to linearized KdV equation} \\ \vspace{0.3cm} 
by 
{\bf \large  Masaki Kawamoto. 
 $^{1}$} and 
 {\bf \large Hisashi Morioka $^{2}$}\\  
Graduate School of Science and Engineering, Ehime University, 3 Bunkyo-cho Matsuyama, Ehime 790-8577. Japan. \\ 
${}^1$ Email: kawamoto.masaki.zs@ehime-u.ac.jp, \\ 
${}^2$ Email: morioka@cs.ehime-u.ac.jp. 
\end{flushleft}

\begin{center}
\begin{minipage}[c]{400pt}
{\bf Abstract} {\small The linearized Korteweg--De Vries equation can be written as a Hamilton-like system. However, the Hamilton energy depends on the time, and is a nonsymmetric operator on $L^2({\bf R})$. By performing suitable unitary transforms on the Hamilton energy, we can reduce this operator into one that is not independent on the time but nonsymmetric. In this study, we consider the $L^2$-stability issues and smoothing estimates for this operator, and prove that it has no eigenvalues.  
}
\end{minipage}
\end{center}

\begin{flushleft}
\end{flushleft}
\begin{flushleft}
\end{flushleft}

\section{Introduction} 
In this study, we consider the $L^2$ properties to obtain a solution to the {\em linearized Korteweg--De Vries (KdV) equation}. Let $x \in {\bf R}$. The linearized KdV operator acting on $L^2({\bf R})$ is defined as follows: 
\begin{align*}
H := H_0 + 12 pV, \qquad H_0 := - p^3 -4 p, 
\end{align*}
where $p = -i \partial _x$, and $V$ is the multiplication operator of $\cosh ^{-2} x$. In a previous work, Pego--Weinstein \cite{PW2, PW} assigned a large scope for $H$ (\cite{PW}, considered more generalized cases including the generalized KdV (gKdV) equation), and showed the absence of embedded eigenvalues for the $H$ and $H^1$ properties of $e^{-itH}$. These results have been applied to KdV and gKdV equations in many studies. On the other hand, one of the important questions in the study of this operator is whether the $L^2$-stable property holds. With regard to the generalized potential case, if the $L^{\infty}$ norm of $V$ and $V'$ are sufficiently small, then this property can be proved using Kato's smooth perturbation method (Kawamoto \cite{Ka}). However, to the best of the authors' knowledge, this issue has not been considered under the ideal potential $V = \cosh ^{-2} x $. Therefore, the purpose of this study is to prove the following theorem: 
\begin{Thm}\label{T1}
For all $t \in {\bf R}$ and $\phi \in L^2({\bf R})$, there exist $0< c_0 <C_0 $ and $0< c_0 ^{\ast} < C_0 ^{\ast}$, which are independent of $t$ such that 
\begin{align*}
c_0 \left\|  \phi \right\| \leq  \left\| e^{-itH} \phi \right\| \leq C_0 \| \phi \|, \quad 
c_0^{\ast} \left\|  \phi \right\| \leq  \left\| e^{-itH^{\ast}} \phi \right\| \leq C_0^{\ast} \| \phi \|
\end{align*}
holds; where $\| \cdot\|$ denotes the norm on $L^2({\bf R})$. 
\end{Thm} 
This theorem enables the consideration of scattering problems for $e^{-itH}$. For $e^{-itH}\phi $, it is important to determine whether the asymptotics 
\begin{align} \label{1}
\lim_{t \to \pm \infty} \left\| 
e^{-itH} \phi - e^{-itH_0} u_{\pm}
\right\| = 0
\end{align}
hold or not, as $pV$ can be regarded as the perturbation of $H_0$. Because it holds that 
\begin{align*}
C_0 \left\| \phi - e^{itH} e^{-itH_0} u_{\pm} \right\| \geq \left\| 
e^{-itH} \phi - e^{-itH_0} u_{\pm}
\right\|, 
\end{align*}
the asymptotics \eqref{1} can be shown by proving the existence of the {\em wave operators} $\CAL{W}^{\pm}$ or their inverse $\CAL{W}^{\pm}_{\mathrm{In}}$
\begin{align*} 
\CAL{W}^{\pm} := \mathrm{s-}\lim_{t \to \pm \infty} e^{itH} e^{-itH_0}, \quad \CAL{W}^{\pm}_{\mathrm{In}} := \mathrm{s-}\lim_{t \to \pm \infty} e^{itH_0} e^{-itH}. 
\end{align*}
In this sense, we consider the existence of wave operators. We obtain the following theorem: 
\begin{Thm}\label{T2}
Wave operators $\CAL{W}^{\pm}$ exist and their inverse $\CAL{W}^{\pm}_{\mathrm{In}}$ also exist.
\end{Thm}   
\begin{Rem}
By virtue of this theorem, we can expect that the asymptotic completeness $\mathrm{Ran}(\CAL{W}^{\pm}) = L^2({\bf R})$ holds. However, we have no idea how to prove $\sigma_{\mathrm{sc}} (H) = \emptyset$. Thus, we omit this discussion, where $\sigma_{\mathrm{sc}} (H)$ indicates a set of singular continuous spectrum of $H$. 
\end{Rem}
The key estimates show that these are the smoothing estimates for $e^{-itH_0}$ and $e^{-itH}$. The smoothing estimates for $e^{-itH_0}$ can be proven as the direct consequence of the results of Ruzhansky--Sugimoto \cite{RS}. Extending these estimates to $e^{-itH} $, we obtain the following theorem:
\begin{Thm}\label{T3}
For all $\alpha >0 $ and $\phi \in L^2({\bf R})$, there exist $C_{1,\alpha}, C_{2, \alpha} >0$ such that 
\begin{align*}
\int_0^{\infty} \left\| e^{- \alpha |x|}  e^{\mp it H} \phi  \right\|^2 dt \leq C_{1,\alpha} \| \phi \|^2
\end{align*}
and 
\begin{align}\label{10}
\int_0^{\infty} \left\| e^{- \alpha |x|} p e^{\mp it H} \phi  \right\|^2 dt \leq C_{2,\alpha} \| \phi \|^2
\end{align}
hold.
\end{Thm}
We also consider the eigenvalue problems for $H$. Let $u: {\bf R} \to {\bf C}$ and $\lambda \in {\bf C}$ such that $H u = \lambda u$ and denote the eigenfunction and eigenvalue of $H$, respectively. Because of the effect of $-4p$ on $H_0$, we can obtain $\sigma (H_0) = \sigma _{\mathrm{ac} }(H_0) = {\bf R}$, where $\sigma _{\mathrm{ac}} (H_0)$ indicates a set of the absolutely continuous spectrum of $H_0$. Clearly, $pV (H_0 + i) ^{-1}$ is a compact operator, and $\sigma _{\mathrm{ess}} (H) = \sigma _{\mathrm{ess}} (H_0)$ holds. Hence, we are interested in the existence of eigenvalues embedded in the real line ${\bf R}$. With regard to previous works \cite{PW2, PW} and references therein, it was shown that $\sigma _{\mathrm{pp}} (H) \cap {\bf R} \subset \{ 0 \} $, which implies that possible eigenvalues that $H$ has is only $0$. However, to the best of the authors' knowledge, there are no results confirming whether $ 0 $ is an eigenvalue or not. Thus, from the theorems in \ref{T1} and \ref{T3}, we present the following theorem and confirm that $H$ has no eigenvalues in ${\bf C}$:
\begin{Thm}\label{T4}
$H$ has no eigenvalues, i.e., if there exist $\lambda \in {\bf C} $ and $u \in L^2({\bf R})$ so that $Hu = \lambda u$, then $u \equiv 0$.
\end{Thm}
In physics, there is a well-known approach that uses the {\em inverse scattering transform} (Leblond \cite{L}). By employing the Jost solutions, the $L^4$-norm of the solutions to 
\begin{align} \label{2}
i \partial_ t u = H u
\end{align}
can be estimated. Recently, studies were conducted on \eqref{2} and mathematical problems, such as Mann \cite{M}, Kato--Kawamoto--Nanbu \cite{KKN},and \cite{Ka}). If one obtains the $L^2$-stable property for \eqref{2}, then we may be able to develop based on these studies.

\section{Auxiliary estimates}
In this section, we introduce the propagation and smoothing estimates for $e^{-itH_0}$ and $e^{-itH}$, which play a key role in analyzing $e^{-itH}$. We denote $\| \cdot \|$ as the norm on $L^2({\bf R})$. First, we introduce the smoothing estimates for $e^{-itH_0}$ (Theorem 5.4 and Corollary 5.5 of Ruzhansky--Sugimoto \cite{RS}).
\begin{Lem}\label{L1}
For all $\phi \in L^2({\bf R})$ and $0 \leq \theta \leq 1$, there exists a constant $C>0$ such that 
\begin{align}\label{3}
\int_{\bf R} \left\| 
\J{x}^{-1} \J{p}^{\theta} e^{\mp itH_0} \phi
\right\| ^2 dt\leq C \| \phi \| ^2 
\end{align}
and 
\begin{align}\label{6}
\int_{\bf R} \left\| 
\J{x}^{-1} |p| e^{\mp itH_0} \phi
\right\| ^2 dt\leq C \| \phi \| ^2 
\end{align}
hold.
\end{Lem}
\Proof{
As \eqref{3} is slightly different from the result of Theorem 5.4 in \cite{RS}, we provide a sketch to explain this. Because 
\begin{align*}
\J{x}^{-1} \J{p}^{\theta} \J{p}^{-1}\J{x}  = \J{x}^{-1} \J{p}^{\theta -1} (x+ i) \cdot (x+ i)^{-1} \J{x},
\end{align*}
and the fact that the operator norm of $(x+i)^{-1} \J{x}$ and 
\begin{align*}
 \J{x}^{-1} \J{p}^{\theta -1} (x+ i) = -i (\theta -1) \J{x}^{-1} p\J{p}^{\theta -3} +   \J{x}^{-1}  (x+ i)  \J{p}^{\theta -1}
\end{align*}
are bounded, we obtain 
\begin{align*}
\int_{\bf R} \left\| 
\J{x}^{-1} \J{p}^{\theta} e^{\mp itH_0} \phi
\right\| ^2 dt\leq C \int_{\bf R} \left\| 
\J{x}^{-1} \J{p}^{} e^{\mp itH_0} \phi
\right\| ^2 dt\leq C \| \phi \|^2, 
\end{align*}
where we use Theorem 5.4 of \cite{RS} with $s=1$ and $a(\xi) = \pm ( \xi ^3 + 4 \xi )$.
}
Next, we introduce an important propagation estimate according to Proposition 4.1 of \cite{PW}, and extend this result to the case of $e^{-itH}$.  
\begin{Lem}\label{L2}
For all $\alpha \in {\bf R}$, $n \in {\bf N} \cup \{ 0 \}$, and $\phi \in \D{e^{\pm \alpha x}}$, 
\begin{align} \label{4}
\left\| p^n e^{\pm \alpha x} e^{\mp it H_0} \phi \right\| \leq 
C t^{-n/2} e^{- \alpha (4-\alpha ^2) t}
\left\| 
e^{ \pm \alpha x} \phi
\right\|
\end{align}
\end{Lem}
\Proof{
Because 
\begin{align*}
e^{\pm \alpha x} e^{\mp it H_0} e^{\mp\alpha x} = e^{\mp it \hat{H}_0}
\end{align*}
with 
\begin{align*}
\hat{H}_0 = -p^3 + 3 \alpha ^2 p -4 p + i \left( \mp 3 p^2 \pm \alpha ^3 \mp 4 \alpha \right),
\end{align*}
we can apply the same argument as in the proof of Proposition 4.1 of \cite{PW}.
}
Because of this lemma, there is the following proposition: 
\begin{Prop}\label{P1}
Let $0 < \alpha \leq 1$. Then, for all $\phi \in L^2({\bf R})$, there exists $C>0$ such that 
\begin{align}\label{7}
\int_{0}^{\infty} \left\| 
e^{-\alpha |x| }  e^{\mp it H} \phi
\right\| ^2 dt \leq C \| \phi \|^2
\end{align}
holds.
\end{Prop}
\Proof{
By Duhamel's formula, one obtains for all $\phi \in L^2({\bf R})$
\begin{align*}
e^{\mp itH} \phi &= e^{\mp itH_0} \phi \mp   12 i \int_0^t e^{\mp i (t-s) H_0} pV e^{ \mp is H} \phi ds,
\end{align*}
which gives 
\begin{align*}
\left\| e^{-\alpha |x| }  e^{\mp itH} \phi \right\| \leq \left\| e^{-\alpha |x| } e^{\mp itH_0} \phi \right\|  + 12
\int_0^t \left\| e^{-\alpha |x|}e^{\mp i (t-s) H_0} pV e^{ \mp is H} \phi \right\|  ds. 
\end{align*}
Employing inequalities
\begin{align*}
& \left\| e^{-\alpha |x|}e^{\mp i (t-s) H_0}, pV e^{ \mp is H} \phi \right\| \\ &
\leq \left\| e^{-\alpha |x| } e^{\mp \alpha x} \right\|_{\SCR{B}} \left\| e^{\pm \alpha x} p  e^{\mp i (t-s) H_0} e^{\mp \alpha x}\right\|_{\SCR{B}}
 \left\| e^{\pm \alpha x} V e^{\alpha |x|} \right\|_{\SCR{B}} \left\| e^{-\alpha |x|} e^{\mp is H} \phi \right\| \\ & \leq 
 C \left( \left\| e^{\pm \alpha x}   e^{\mp i (t-s) H_0} e^{\mp \alpha x}\right\|_{\SCR{B}} + \left\| p e^{\pm \alpha x}   e^{\mp i (t-s) H_0} e^{\mp \alpha x}\right\|_{\SCR{B}} \right)\left\| e^{-\alpha |x|} e^{\mp is H} \phi \right\| 
\end{align*}
and using Lemma \ref{L2} with $n=0$ and $n=1$, we obtain 
\begin{align} \nn 
& \left\| e^{-\alpha |x| }e^{\mp itH} \phi \right\| \\ & \leq \left\| e^{-\alpha |x| }e^{\mp itH_0} \phi \right\|  + C,
\int_0^t  (1 + (t-s)^{-1/2}) e^{-\alpha (4-\alpha ^2) (t-s)}   \left\| e^{- \alpha |x|} e^{ \mp is H} \phi \right\| ds, \label{5}
\end{align}
where $\| \cdot \|_{\SCR{B}}$ denotes the operator norm on $L^2({\bf R})$. Because there exists a $t$-independent constant $C >0$ such that 
\begin{align*}
\int_0^t (t-s)^{-1/2} e^{-\alpha (4-\alpha ^2) (t-s)} ds \leq \int_0^1 \tau ^{-1/2} d \tau + \int_1^{t} e^{-\alpha (4- \alpha ^2) \tau} d \tau \leq C
\end{align*}
holds, Glonwall's inequality and \eqref{5} yield a $t$-independent constant $C>0$ such that
\begin{align} \label{9}
\left\| 
 e^{-\alpha |x| }e^{\mp itH} \phi 
\right\| \leq \left\|  e^{-\alpha |x| }e^{\mp itH_0} \phi  \right\| +  C \int_0^t \left\|  e^{-\alpha |x| }e^{\mp isH_0} \phi  \right\| I_1(t,s) ds, 
\end{align}
where 
\begin{align*}
I_1 (t,s) &:= (1 + (t-s) ^{-1/2}) e^{-\alpha (4-\alpha ^2) (t-s)} \mathrm{exp} \left( \int_s^t (1 + (t- \tau)^{-1/2}) e^{-\alpha (4-\alpha ^2) (t-\tau)} d \tau  \right) \\ 
& \leq C (1 + (t - s ^{-1/2}) e^{-\alpha (4-\alpha ^2) (t-s)}.
\end{align*}
We note that 
\begin{align*}
 \int_s^t (1 + (t- \tau)^{-1/2}) e^{-\alpha (4-\alpha ^2) (t-\tau)} d \tau \leq  \int_0^t (1 + (t- \tau)^{-1/2}) e^{-\alpha (4-\alpha ^2) (t-\tau)} d \tau \leq C.
\end{align*}
Inequality \eqref{9} yields 
\begin{align*}
\left\| 
 e^{-\alpha |x| }e^{\mp itH} \phi 
\right\|^2 \leq C \left\|  e^{-\alpha |x| }e^{\mp itH_0} \phi  \right\|^2 + C \left( \int_0^t \left\|  e^{-\alpha |x| }e^{\mp isH_0} \phi  \right\| \left( I_2(t,s) + I_3(t,s) \right) ds \right)^2
\end{align*}
with 
\begin{align*}
I_2(t,s) := e^{-\alpha (4-\alpha ^2) (t-s)}, \quad I_3(t,s) := (t-s)^{-1/2} e^{-\alpha (4-\alpha ^2) (t-s)}. 
\end{align*}
To show this proposition, it suffices to show 
\begin{align*}
\int_0^{\infty} \left( \int_0^t \left\|  e^{-\alpha |x| }e^{\mp isH_0} \phi  \right\| \left( I_2(t,s) + I_3(t,s) \right) ds \right)^2 dt \leq C \|\phi \|^2. 
\end{align*} 
 We now estimate the term associated with $I_3(t,s)$ because $I_2(t,s)$ can be estimated by the same approach. With the Cauchy--Schwarz inequality,
 we have 
 \begin{align*}
& \int_0^{\infty} \left( \int_0^t \left\|  e^{-\alpha |x| }e^{\mp isH_0} \phi  \right\|I_3(t,s) ds \right)^2 dt  \\ & \leq
 C \int_0^{\infty} \left( \int_0^t \left\|  e^{-\alpha |x| }e^{\mp isH_0} \phi  \right\|^2 I_3(t,s) ds  \right) \left( \int_0^t I_3(t, \tau) d \tau\right) dt 
 \\ & \leq 
 C \int_0^{\infty} \int_0^t \left\|  e^{-\alpha |x| }e^{\mp isH_0} \phi  \right\|^2 I_3(t,s) ds dt \\ 
 &= C \int_0^{\infty} \int_s^{\infty}  \left\|  e^{-\alpha |x| }e^{\mp isH_0} \phi  \right\|^2 I_3(t,s) dt ds \\ 
 &= C \int_0^{\infty}  \left\|  e^{-\alpha |x| }e^{\mp isH_0} \phi  \right\|^2 \left( \int_s^{\infty}   I_3(t,s) dt \right) ds  \\ & \leq C \| \phi \|^2; 
 \end{align*}
hence, the proof is completed, where we use
\begin{align*}
\int_s^{\infty} I_3 (t, s) d t &= \int_s^{s+1} I_3 (t, s) d t + \int_{s+1}^{\infty} I_3(t, s) d t \\ &\leq \int_s^{s+1} (t-s)^{-1/2} dt + \int_{s+1}^{\infty} e^{- \alpha (4-\alpha ^2) (t-s)} dt \leq C.
\end{align*}
}
\begin{Rem}
After the commutator calculations, we can obtain the smoothing effect
\begin{align*}
\int_{\bf R} \left\| 
e^{-\alpha |x| } p  e^{\mp it H} \phi
\right\| ^2 dt \leq C \| \phi \|^2. 
\end{align*}
However, for simplicity, we omitted the proof, which is provided in \S{4}; indeed, this property is not used in the proof of our main theorems. 
\end{Rem}

\section{Proof of Theorems \ref{T1}, \ref{T2}, and \ref{T4}}
We now prove the main theorems. First, we introduce the following lemma: 
\begin{Lem}\label{L3}
For all $t \in {\bf R}$ and $\phi \in L^2({\bf R})$, 
\begin{align*}
\left\| 
e^{-itH} \phi
\right\| \leq e^{5 |t|} \| \phi \|
\end{align*}
holds.
\end{Lem}
\Proof{
Decompose $H$ into 
\begin{align*}
H = \tilde{H}_0 -6iV',
\end{align*}
where
\begin{align*}
\tilde{H}_0 := -p^3 - 4p + 6(pV + Vp)
\end{align*}
and $V'$ is the multiplication operator of $V'(x) = -2 \sinh x /\cosh ^3 x$. Clearly, $(pV + Vp) (-p^3 -4p + i)^{-1}$ is the compact operator. We notice that $\tilde{H}_0$ is the self-adjoint operator on $\D{\tilde{H}_0} = H^3({\bf R})$. Then, the Duhamel formula gives 
\begin{align*}
e^{-itH} \phi = e^{-it \tilde{H}_0} \phi - 6\int_0^t e^{-i(t-s) \tilde{H}_0} V' e^{-isH} \phi ds,
\end{align*}
which holds for all $\phi \in L^2({\bf R})$. This equation, unitarity of $e^{-it \tilde{H}_0}$ on $L^2({\bf R})$, $6 \sup |V'(x)| \leq 5$, and Gronwall's inequality imply 
\begin{align*}
\left\| 
e^{-itH} \phi
\right\| \leq e^{5 |t| } \| \phi \|.
\end{align*}
}
\subsection{Proof of Theorem \ref{T1}}
Let $u,v \in \SCR{S} ({\bf R})$, $0< \alpha \leq 1$, and $(\cdot , \cdot)$ denote the inner product on $L^2({\bf R})$. Because 
\begin{align*}
\frac{d}{dt} \left( 
e^{\pm itH_0} e^{\mp itH} u , v
\right) =  \mp  12i \left( e^{\pm itH_0} pV e^{\mp itH} u,v \right), 
\end{align*}
we have 
\begin{align}
\nn & \left| \left( 
e^{\pm itH_0} e^{\mp itH} u , v
\right) \right| \\ 
\nn &\leq \| u \| \left\|  v \right\|  +  12 \int_0^{t}\left|  \left(  pV e^{\mp isH} u, e^{\mp isH_0}v \right) \right| ds 
\\ 
\nn & \leq  \| u \| \| v \| + C \int_0^t \left\| e^{-\alpha |x| } e^{\mp is H} u \right\| \left\| e^{-\alpha |x|} p e^{\mp is H_0} v  \right\| ds \\ 
\nn & 
\leq  \| u \| \| v \| + C \left( \int_{\bf R} \left\| e^{-\alpha |x| } e^{\mp is H} u \right\|^2 dt \right) ^{1/2} \left(  \int_{\bf R} \left\| e^{-\alpha |x|} p e^{\mp is H_0} v  \right\| ^2 ds \right)^{1/2} \\ 
\label{8} & \leq C \| u \| \|  v\|, 
\end{align}
where we use the Cauchy--Schwarz inequality, \eqref{6} and \eqref{7}. Hence, we obtain $u \in \SCR{S} ({\bf R})$
\begin{align*}
\left\| 
e^{\mp itH} u
\right\| &= \left\| e^{\pm itH_0} e^{\mp it H} u \right\|
= \sup_{\| v \| = 1} \left| \left( 
e^{\pm itH_0} e^{\mp itH} u , v
\right) \right|  \leq C \| u \|. 
\end{align*}
Owing to Lemma \ref{L3} and the density argument, we can extend this estimate to all $\phi \in L^2({\bf R})$ and show 
that for all $t \in {\bf R}$ and $\phi \in L^2({\bf R})$, there exists a $t$-independent constant $C_0>0$ such that 
\begin{align*}
\left\| 
e^{ -itH} \phi
\right\| \leq C_0 \| \phi \|
\end{align*}
holds. Conversely, 
\begin{align*}
\left\| \phi \right\| = \left\| e^{itH} e^{-itH} \phi \right\| \leq C_0 \left\| e^{-itH} \phi \right\|
\end{align*}
implies $c_0 \| \phi \| :=(1/C_0) \| \phi \| \leq \| e^{-itH} \phi \| $. 

As for $e^{-itH^{\ast}}$, we use 
\begin{align*}
e^{\mp itH^{\ast}} u &= e^{\mp itH_0} u \mp  12 i \int_0^t e^{\mp i(t-s) H_0} Vp e^{\mp is H^{\ast}} u ds, \\ 
&= e^{\mp itH_0} u \mp 12 i \int_0^t e^{\mp i(t-s) H_0} pV e^{\mp is H^{\ast}} u ds  \pm  12 \int_0^t e^{\mp i(t-s) H_0} V' e^{\mp is H^{\ast}} u ds,
\end{align*}
which enables us to obtain the estimate 
\begin{align*}
\int_0^{\infty} \left\| e^{-\alpha |x|} e^{\mp itH^{\ast}} u \right\|^2 dt \leq C \| u \|^2.
\end{align*}

\subsection{Proof of Theorem \ref{T2}}
We present proof of Theorem \ref{T2}. The existence of 
\begin{align*}
\mathrm{s-}\lim_{t \to \pm \infty} e^{itH_0}e^{-itH}
\end{align*}
can be proven by \eqref{8}. Hence, we only prove the existence of $\CAL{W}^{\pm}$. Owing to Theorem \ref{T1}, we have $e^{itH}e^{-itH_0}$, which is uniformly bounded in $t$. According to the density argument, it suffices to show the existence of $\lim_{t \to \pm \infty} e^{-itH}e^{-itH_0} \phi$, $\phi \in \SCR{S} ({\bf R})$, which is immediately proven by 
\begin{align*}
\left\| \frac{d}{dt} e^{itH} e^{-itH_0} \phi \right\| \leq 12 C_0 \left\| pV e^{-itH_0} \phi \right\| \leq C (1 + t^{-1/2} )e^{-\alpha (4-\alpha ^2)t } \left\| e^{\alpha x} \phi \right\|
\end{align*}
and the Cook--Kuroda method.

\subsection{Proof of Theorem \ref{T4}}
We now prove Theorem \ref{T4}, which is quite simple. First, we show $\sigma_{\mathrm{pp}} (H) \cap \left( {\bf C} \backslash {\bf R}  \right) = \emptyset$ Suppose $\lambda = a + i b$, $a \in {\bf R}$, and $b \in {\bf R} \backslash \{ 0\}$ as the eigenvalue of $H$, and $u$ as the associated eigenfunction with $\| u \| \neq 0$. Then, $\left\| e^{-itH} u \right\| = e^{bt} \| u \|$ tends to infinity or $0$ as $t \to \infty$, which contradicts Theorem \ref{T1} and implies that $u \equiv 0$. Second, we show $\sigma_{\mathrm{pp}} (H) \cap {\bf R}  = \emptyset$. Because Proposition \ref{P1} holds {\em for all $\phi \in L^2({\bf R})$}, we take $\phi \in L^2({\bf R})$ such that $H \phi = \lambda \phi$ with $\lambda \in {\bf R}$. Then, \eqref{7} implies  
\begin{align*}
 \int_0^{\infty} \left\| e^{- \alpha |x|} \phi \right\| ^2 dt = \int_0^{\infty} \left\| e^{- \alpha |x|} e^{-it H} \phi \right\| ^2 dt \leq C \| \phi \|^2.
\end{align*}
This inequality is true only in the case where $\phi  \equiv 0$; otherwise, it false.

\section{Smoothing estimate for $e^{-itH}$}
In this section, we present the smoothing estimate \eqref{10}. Let $0 < \alpha \leq 1$. By the Duhamel formula, we have 
\begin{align*}
& e^{-\alpha |x|} (p+ 2i) e^{-itH} u \\  &= e^{-\alpha |x|} (p+ 2i) e^{-itH_0} u \\ & \quad  + 12 \int_0^t e^{- \alpha |x|} e^{-i(t-s) H_0} p e^{-\alpha |x|} \cdot e^{\alpha |x|} (p+ 2i) V (p+ 2i)^{-1} e^{\alpha |x|} \cdot e^{-\alpha |x|} (p+ 2i) e^{- is H} u ds. 
\end{align*}
Hence, we now estimate the operator norm of $J_1 := e^{\alpha |x|} (p+ 2i) V (p+ 2i)^{-1} e^{\alpha |x|}$. Clearly, $J_1 = J_2 -i J_3$ holds with $J_2 := e^{2\alpha |x| } V $ and $J_3 = e^{\alpha |x| } V'(x) (p+ 2i)^{-1} e^{\alpha |x|} $. Because $\| J_2 \|_{\SCR{B}} \leq C $ and $\| e^{\alpha |x| } V'(x) \cosh x  \|_{\SCR{B}}  \leq C$ holds, we estimate 
\begin{align*}
\left\| 
J_1 
\right\|_{\SCR{B}} & \leq C  + \left\| \left(  \cosh^{-1} x \right) (p+ 2 i)^{-1} e^{\alpha |x|}\right\|_{\SCR{B}} 
\\ & =   C  + \left\|  e^{\alpha |x|} (p - 2 i)^{-1}  \left( \cosh^{-1} x \right) \right\|_{\SCR{B}}
\\ & \leq 
C  + \left\|  e^{\alpha x} (p - 2 i)^{-1}  \left( \cosh^{-1} x \right) \right\|_{\SCR{B}} + \left\|  e^{-\alpha x} (p - 2 i)^{-1}  \left( \cosh^{-1} x \right) \right\|_{\SCR{B}} \\ & \leq 
C + C \left\| e^{\pm \alpha x} (p \mp 2i )^{-1} e^{\mp \alpha x} \right\|_{\SCR{B}} \left\| e^{\alpha |x|}  \left( \cosh^{-1} x \right)  \right\|_{\SCR{B}}.
\\ & \leq  
C + C \left\|  (p \mp ( 2 - \alpha) i )^{-1} \right\|_{\SCR{B}}. 
\\ & \leq C, 
\end{align*}
where we use $0 < \alpha \leq 1$, i.e., $|2- \alpha| \geq 1$. By $\| J_1\|_{\SCR{B}} \leq C$ and the same argument in proving Proposition \ref{P1}, we have 
\begin{align*}
\left\| 
 e^{-\alpha |x| } (p+ 2i)e^{\mp itH} \phi 
\right\|^2 &\leq C \left\|  e^{-\alpha |x| } (p+ 2i)e^{\mp itH_0} \phi  \right\|^2 \\ & \quad + C \left( \int_0^t \left\|  e^{-\alpha |x| }(p+ 2i)e ^{\mp isH_0} \phi  \right\| \left( I_2(t,s) + I_3(t,s) \right) ds \right)^2. 
\end{align*}
This yields the smoothing 
\begin{align*}
\int_0^{\infty} \left\| 
 e^{-\alpha |x| } (p+ 2i)e^{\mp itH} \phi 
\right\|^2 dt \leq C \| \phi  \|^2. 
\end{align*}
We note that $e^{-\alpha |x| } p = e^{- \alpha |x|} p(p+2i)^{-1} e^{\alpha |x|} \cdot e^{- \alpha |x|} (p+ 2i)$, and that
\begin{align*}
\left\|
e^{- \alpha |x|} p(p+2i)^{-1} e^{\alpha |x|}
 \right\|_{\SCR{B}} \leq C 
\end{align*}
holds by using a similar argument to show that $\|  J_1\|_{\SCR{B}} \leq C $, \eqref{10} can be proven. 

\end{document}